\documentclass[a4paper]{article}

\author{Joaquim Ro\'e \\
\small{Departament de Matem\`{a}tiques, Universitat Auton\`{o}ma de
Barcelona,}\\ \small{Edifici C, 08193-Bellaterra
  (Barcelona), Spain. e-mail: jroe@mat.uab.es}}

\usepackage[centertags]{amsmath}
\usepackage{amsthm}
\usepackage{amssymb}

\theoremstyle{definition}
\newtheorem{Def}{Definition}
\theoremstyle{plain}

\newtheorem{Cor}[Def]{Corollary}

\newtheorem{Teo}[Def]{Theorem}

\theoremstyle{remark}

\newcommand{\Z}{\mathbb{Z}}

\renewcommand{\P}{\mathbb{P}}

\begin{document}

\title{On submaximal plane curves}
\maketitle

\begin{abstract}
We prove that a submaximal curve in $\P^2$ has sequence of
multiplicities $(\mu,\nu, \dots, \nu)$, with $\mu < s \nu $ for
every integer $s$ with $(s-1)^2(s+2)^2\ge 6.76(\, r-1)$.
\end{abstract}


\bigskip

This note is a sequel to \cite{Roe01}, where a specialization
method was developed in order to bound the degree of singular
plane curves. The problem under consideration is, given a system
of multiplicities $(m)=(m_1,m_2,\dots,m_r) \in \Z^r$ and points
$p_1, \dots, p_r \in \P^2$, which we shall always assume to be in
general position, to determine the minimal degree $\alpha(m)$ of a
curve with multiplicity $m_i$ at each point $p_i$. In
\cite{Roe01}, the focus was on homogeneous $(m)$, (i.e.,
$m_1=m_2=\dots=m_r$), but the method applies in general; here it
is used to show that if one of the multiplicities is much bigger
than the others, in a sense we make precise below (see theorem
\ref{nohomo}), then
\begin{equation}
\label{eqnag} \alpha(m)> \frac{\sum_{i=1}^r m_i}{\sqrt{r}}.
\end{equation}

In connection with his solution to the fourteenth problem of
Hilbert, M. Nagata conjectured in 1959 that the inequality
(\ref{eqnag}) holds for all $(m)$ provided $r>9$, and proved it in
the case when $r$ is a perfect square (see \cite{Nag59} or
\cite{Nag60}). Since then, many partial results have been proved
by several authors (see for instance \cite{Xu94}, \cite{Eva98},
\cite{Har01}, \cite{HR??}, \cite{Sze01}, \cite{Tut??}), but as far
as we know the conjecture remains open in general. One of the
research lines in this area is the study of \emph{submaximal
curves} that arise in the context of Seshadri constants. A
submaximal curve is an irreducible curve which causes the
($r$-point) Seshadri constant of a surface to be non-maximal; in
the case of $\P^2$ it is just an irreducible curve which causes
(\ref{eqnag}) to fail, and Nagata's conjecture can be equivalently
stated by saying that there exist no submaximal curves for $r>9$.
T. Szemberg proved in \cite[4.6]{Sze01} that every submaximal
curve on a surface with Picard number $\rho=1$ whose multiple
points are in general position is quasi-homogeneous, i.e., has
$m_2=\dots=m_r$ for a suitable ordering of the points. In the case
of the projective plane, our result shows that moreover $m_1$ can
not be much bigger than $m_2$, constraining further the range of
possible counterexamples to Nagata's conjecture. It is worth
mentioning that quasi-homogeneous curves are relevant also for the
method of C. Ciliberto and R. Miranda \cite{CM98a} to compute the
dimension of (homogeneous) linear systems.



\medskip

The approach is based on the specialization introduced in
\cite{Roe01}. Roughly speaking, one proves that if there exists a
curve with given multiplicities at $r$ general points, then by
semicontinuity there must also exist curves with the same degree
and (virtual) multiplicities at $r$ points which satisfy some well
chosen proximity relations. The proximity inequalities impose then
that the effective multiplicity of the specialized curve must
grow, and one uses this bigger multiplicity as a bound for the
degree of the curve.

In order to give a brief explicit description of the
specialization, let us recall the notations of \cite{Roe01} (see
\cite{Cas00} for generalities on clusters and unloading, and
\cite{Roe?4} for a general approach to specializations
parameterized by varieties of clusters). We work on $\P^2$, and
consider both proper and infinitely near points (which are those
lying on a smooth surface that dominates $\P^2$ birrationally. A
\emph{cluster} is a set $K$ of points of $\P^2$ such that if $p\in
K$ and $p$ is infinitely near to $q$ (i.e., $p$ lies on the
exceptional divisor of $q$ after blowing up a sequence of points)
then $q \in K$. Write $\pi_K:S_K\longrightarrow \P^2$ for the
blowing up of all points in $K$. For each $i=2, \dots, r$, we
denote $U_{i}$ the set of clusters $K=\{p_1, \dots, p_r \}$ such
that
\begin{itemize}
\item $p_2, \dots, p_i$ are proximate to $p_1$,
\item $p_j$ is proximate to $p_{j-1}$ for all $j=2, \dots, r$, and
\item there are no other proximity relations.
\end{itemize}
In other words, denoting by $E_j$ the (total) exceptional divisor
of blowing up $p_j$ this can be expressed by saying that the
divisors $\tilde E_1= E_1-E_2-\dots-E_i$ and $\tilde
E_j=E_j-E_{j+1}$, $j=2, \dots r-1$ on the surface $S_K$ are
effective and irreducible.

The sets $U_{i}$ are nonempty and have a natural structure of
smooth irreducible locally closed subvarieties in a projective
variety (the iterated blowing-up $X_{r-1}$ of Kleiman
\cite{Kle81a}), and they satisfy
$$\overline{U_{2}} \supset \overline{U_{3}}\supset \dots \supset \overline{U_{r}}.$$
The specialization we use works stepwise. Begin with a general
cluster of $r$ distinct points, $K=\{p_1,p_2,\dots,p_r\}$, and a
curve $C$ with multiplicities $(m)=(m_1, m_2, \dots, m_r)$ at
these points, assuming $m_1 \ge m_2 \ge \dots \ge m_r$. Then
specialize $K$ to a cluster $K_3$ general in $U_{3}$. If
$m_1<m_2+m_3$, then the specialized curve $C_3$ cuts negatively
the irreducible divisor $\tilde E_1=E_1-E_2-E_3$, so $\tilde E_1$
is a component of $\pi_K^*(C_3)-m_1E_1-\dots-m_rE_r$, and the
effective multiplicity of $C_3$ at $p_1$ is bigger than $m_1$.
Call $(m^{(3)})$ the system of multiplicities obtained after
unloading multiplicities (i.e., substracting the $\tilde E_j$
which are cut negatively); $C_3$ goes through the cluster $K_3$
with multiplicities $(m^{(3)})$. Then specialize $K_3$ to a
general $K_4 \in U_4$, and successively to a $K_5 \in U_5$, \dots,
to a $K_r\in U_r$, performing unloadings whenever necessary. The
first multiplicity of the last system $(m^{(r)})$ is a lower bound
for the degree of $C$, and therefore $\alpha(m)\ge m_1^{(r)}$.
This multiplicity is not hard to compute in each particular case;
in \cite{Roe01} a bound was given that holds in general and is
asymptotically sharp but that in many particular cases can be
improved, especially when the multiplicities are relatively small.
Now we are interested in the case that $m_1$ is much bigger than
the other multiplicities, in which one can show that the
inequality (\ref{eqnag}) holds:

\begin{Teo}
\label{nohomo}
  Let $s$ be such that $(s-1)^2(s+2)^2 \ge 6.76(r-1)$,
$r >9$, and assume $m_1 \ge m_2 \ge \dots \ge m_r$. If moreover
$$
m_1 \ge \sum_{i=2}^{s+1} m_i
$$
then $\alpha(m) > \sum_{i=1}^r m_i/\smash{\sqrt{r}}$.
\end{Teo}

\begin{proof}
Using notations as above, the preceding discussion shows that it
suffices to prove $m_1^{(r)} > \sum_{i=1}^r m_i/\smash{\sqrt{r}}$.
The hypothesis implies that the system of multiplicities $(m)$ is
consistent for all clusters in $U_3$, \dots, $U_{s+1}$ (no
unloading is needed for these) so
$(m)=(m^{(3)})=\dots=(m^{(s+1)})$. Then, apply \cite[lemma
3.5]{Roe01} as in the proof of \cite[theorem 4.1]{Roe01} to obtain
$$
 m_1^{(r)} \ge \sum_{i=1}^r m_i \, \left(1-\frac{1}{r}\right)
   \prod_{k = s+1}^{r-1} \left(1 - {k\over {{k^2} + r-1}} \right).
$$
We have to see that this is bigger than $\sum m_i /
\smash{\sqrt{r}}$. Because of \cite[proposition 5.1]{Roe01}, it
will be enough to prove
\begin{equation}
\label{untros}
  \prod_{k = 2}^{s} \left(1 - {k\over {{k^2} + r-1}} \right)^{-1}
  > \frac{\sqrt{r}}{\sqrt{r-1}-\frac{\pi}{8}}.
\end{equation}
Write $x^2=r-1$. As $r >9$, we have $x \ge 3$. The term on the
left in (\ref{untros}) is
$$
 \prod_{k = 2}^{s} \left(1 + \frac {k}{k(k+1) + x^2} \right)>
 1+\sum_{k = 2}^{s} \frac {k}{k(k+1) + x^2}.
$$
Let $s_0$ be the minimum integer such that $(s_0-1)^2(s_0+2)^2 \ge
6.76x^2$; the hypothesis on $s$ says that $s \ge s_0$, and it is
clear that $s_0(s_0+1) \le x^2$. Using this we get
$$
1+\sum_{k = 2}^{s} \frac {k}{k(k+\!1) + x^2}>
 1+\sum_{k = 2}^{s_0} \frac {k}{2 x^2} =
 1+ \frac{(s_0-\!1)(s_0+\!2)}{4 x^2} \ge
 1+ \frac{\sqrt{6.76}x}{4 x^2}= 1 + \frac{.65}{x}.
$$
On the other hand, the term on the right in \ref{untros}
can be written as
$$
\frac{x}{x-\pi/8} \sqrt{1+\frac{1}{x^2}}\le
\left(1+\frac{\pi}{8x-\pi}\right)\left(1+\frac{1}{2x^2}\right)
$$
which for $x \ge 3$ is less or equal to $1+.65/x$, and the proof
is complete.

\end{proof}

\begin{Cor}
\label{quasi} Let $C$ be a submaximal curve with respect to
general points $p_1, p_2 \dots, p_r \in \P^2$. Then for some
reordering of the points, the system $(m)$ of multiplicities of
$C$ at $p_1, p_2 \dots, p_r$ is $(m)=(\mu,\nu, \dots,\nu)$ with
$\mu < s \nu $ for every integer $s$ with $(s-1)^2(s+2)^2\ge
6.76(r-1)$.
\end{Cor}
\begin{proof}
The system $(m)$ of multiplicities of $C$ at $p_1, p_2 \dots, p_r$
is $(m)=(\mu,\nu, \dots,\nu)$ because of \cite[corollary
4.6]{Sze01}, so it is enough to prove that $\mu \ge s \nu $ for
some integer $s$ with $(s-1)^2(s+2)^2\ge 6.76(r-1)$ leads to
contradiction. But theorem \ref{nohomo} shows that there are no
submaximal curves when $\mu \ge s \nu $ for some integer $s$ with
$(s-1)^2(s+2)^2\ge 6.76(r-1)$, so we are done.
\end{proof}

We finish by an example, considering the smallest values $r$ for
which Nagata's conjecture is unknown. Corollary \ref{quasi} says
that a submaximal curve with respect to $r$ general points, $10
\le r \le 15$ has system of multiplicities $(m)=(\mu,\nu,
\dots,\nu)$ with $\mu < 3 \nu $.

\bibliographystyle{amsplain}
\bibliography{Biblio}

\end{document}